\documentclass[a4paper]{article}
\usepackage{amsmath,amssymb,amsthm}
\usepackage[latin1]{inputenc}
\usepackage[all]{xy}
\usepackage[dvips]{graphicx}
\usepackage{fancyhdr}

\newcommand{\NN}{\mathbb{N}}
\newcommand{\ZZ}{\mathbb{Z}}

\newcommand{\kk}{\mathbb{K}}
\newcommand{\rk}{\operatorname{rk}}
\newcommand{\MPi}{\mathrm{M}\Pi}
\newcommand{\Mchi}{\mathrm{M}\chi}
\newcommand{\lbr}{\langle}
\newcommand{\rbr}{\rangle}
\newcommand{\Pe}{\mathcal{P}erm}
\newcommand{\Sy}{\mathbb{S}}
\newcommand{\Pli}{\mathcal{P}re\mathcal{L}ie}
\newcommand{\Lie}{\mathcal{L}ie}
\newcommand{\Mag}{\mathcal{M}ag}
\newcommand{\Comtri}{\mathcal{C}om\mathcal{T}rias}
\newcommand{\Postlie}{\mathcal{P}ost\mathcal{L}ie}
\newcommand{\zero}{\hat{0}}
\newcommand{\one}{\hat{1}}
\newcommand{\pted}[1]{\mathbf{\overline{#1}}}
\newcommand{\FF}{\mathcal{F}}

\newtheorem{theorem}{Theorem}[section]
\newtheorem{proposition}[theorem]{Proposition}

\newtheorem{lemma}[theorem]{Lemma}

\theoremstyle{definition}
\newtheorem{definition}{Definition}

\theoremstyle{remark}
\newtheorem{remark}[theorem]{Remark}

\title{Pointed and multi-pointed partitions \\ of type $A$ and $B$}
\author{F. Chapoton and B. Vallette}
\date{\today}

\begin{document}

\maketitle

\begin{abstract}
  The aim of this paper is to define and study pointed and
  multi-pointed partition posets of type $A$ and $B$ (in the
  classification of Coxeter groups). We compute their characteristic
  polynomials, incidence Hopf algebras and homology groups. As a
  corollary, we get that some operads are Koszul over $\mathbb{Z}$.
\end{abstract}

\section*{Introduction}

For every finite Weyl group $W$, there exists a generalized partition
poset (\emph{cf.} \cite{BW2}) defined through the hyperplane
arrangement of type $W$. In the case $A_{n-1}$, this poset is the
usual poset of partitions of $\{1,\ldots ,\, n\}$. B.  Fresse has
shown in \cite{Fresse} that this poset also arises from the theory of
operads. A pointed and a multi-pointed variation of this poset have
been defined by the second author in \cite{homologie}, once again in
the context of Koszul duality of operads. In this article, we study
the main properties of these two types of posets. One motivation for
this article was the idea that there should also exist a pointed
partition poset and a multi-pointed partition poset for other Weyl
groups. Here we propose a definition for the pointed partition posets
of type $B$ and check that it satisfies most of the properties which
are expected in general and hold in type $A$. This definition has been
guessed by similarity, but we hope that
there is a general definition of geometric nature, to be found.\\

Let us summarize briefly what properties the generalized pointed
partition poset associated to a Weyl group should have. Let $h$ be the
Coxeter number and $n$ be the rank of the Weyl group $W$.  Then its
characteristic polynomial should be $(x-h)^n$; the number of maximal
elements should be $h$, with a transitive action of the Weyl group.
Also the characteristic polynomial of any maximal interval should be
$(x-1)(x-h)^{n-1}$ and the homology must be concentrated in maximal
dimension. We prove that all these properties hold in type $A$ and
$B$.

One can remark that the expected characteristic polynomial is the
same as the characteristic polynomial of the hyperplane
arrangement called the Shi arrangement \cite{athanase,headley}.
One difference is that there is no action of the Weyl group on the
Shi arrangement. Going to the limit where parallel hyperplanes
come together gives the so-called double Coxeter arrangement
\cite{solotera}, which is no longer a hyperplane arrangement in
the usual sense. Still the double Coxeter arrangement is free, and
all its degrees are the Coxeter number. Maybe the pointed
partition poset should be thought of as the missing intersection
poset for the double Coxeter arrangement.

It may be worth noting that there is a family of posets which are
probably related to the pointed partition posets of type $A$, in a
rather non-evident way. It is made of some posets on forests of
labeled rooted trees, introduced by J. Pitman in \cite{pitman}, which
seem to share the same characteristic polynomials. Maybe there is an
homotopy equivalence between these posets of forests and the posets of
pointed partitions.\\

Let us now state what are the expected properties of the
multi-pointed partition posets associated to a Weyl group. Let
$e_1,\dots,e_n$ be the exponents of $W$. Then the characteristic
polynomial should be $\prod_{i=1}^{n}(x-(h+e_i))$ and the homology
of the poset must be concentrated in maximal dimension. This is
true in type $A$. We have been unable to guess what should be the
poset of type $B$, even if we have some evidence that it should
exist.

The multi-pointed partition poset should be related to the
so-called Catalan arrangement \cite{athanase}. The characteristic
polynomials coincide, but the posets are different in general. As
the Weyl group acts on the Catalan arrangement, one may wonder if
the action is the same on the top homology of both posets, which
should have the same dimension. If this is true, this may come
from an homotopy equivalence between the (realizations of the) two
posets.

Just as in the Shi case, one can take the limit of the Catalan
arrangement where parallel hyperplanes come together. This gives
the triple Coxeter arrangement \cite{terao}. This triple Coxeter
arrangement is free, and its degrees are the roots of the
characteristic polynomial of the Catalan arrangement.\\

In type $A$, just as for the usual partition lattices, the pointed and
multi-pointed partitions posets give rise to interesting actions of
the symmetric groups $\Sy_n$ on their homology groups. We use the
relations with the theory of Koszul operads (based on representations
of $\Sy_n$), described by the second author in \cite{homologie}, to
compute this action on the homology groups of the pointed and
multi-pointed posets. In the other direction, the fact that the
various posets are Cohen-Macaulay over $\mathbb{Z}$ implies that some
operads are Koszul over $\mathbb{Z}$, which is an important result in
the study of the deformations of algebraic structures. \\

Since the posets studied here have nice properties for their intervals
and products, we can associate to them an incidence Hopf
algebra.\\

\subsection*{Conventions}

All posets are implicitly finite. A poset $\Pi$ is said to be
\emph{bounded} if it admits one minimal and one maximal element
(denoted by $\zero$ and $\one$). It is \emph{pure} if for any
$x\leqslant y$, all maximal chains between $x$ and $y$ have the
same length. If a poset is both bounded and pure, it is called a
\emph{graded} poset. A pure poset with a minimal element is called
\textit{ranked}.

For a graded or ranked poset $\Pi$ with rank function $\rk$, the
\emph{characteristic polynomial} is defined by the following
formula:
\begin{equation}
  \chi(x):=\sum_{a \in \Pi} \mu(a)x^{n-\rk(a)},
\end{equation}
where $\mu$ denotes the M\"obius function of the poset $\Pi$ and
$n$
is the rank of the maximal elements.\\

Let $\kk$ be the ring $\ZZ$ or any field. Denote by $[n]$ the set
$\{1,\dots,n\}$.

\section{Pointed partition posets}

In this section, we give the definitions of the pointed partition
posets and some of their basic properties.

\subsection{Definitions}

\subsubsection{Type $A$}

\label{definitionPA}

First, we recall the definition, introduced in \cite{homologie},
of the pointed partition poset of type $A_{n-1}$.

\begin{definition}[Pointed partition]
A \emph{pointed partition} of $[n]$ is a partition of $[n]$
together with the choice of one element inside each block, called
the \emph{pointed element} of this block.
\end{definition}

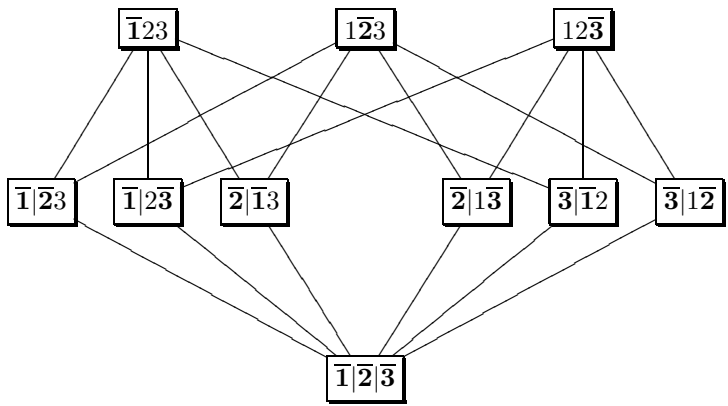
\begin{figure}[h]
$$\xymatrix@R=50pt@C=15pt{ & *+[F-,]{\pted{1}23} \ar@{-}[dl]
\ar@{-}[d] \ar@{-}[dr] \ar@{-}[drrrr] & &*+[F-,]{1\pted{2}3}
\ar@{-}[dl] \ar@{-}[dlll] \ar@{-}[dr] \ar@{-}[drrr] & &
*+[F-,]{12\pted{3}} \ar@{-}[dllll] \ar@{-}[dl] \ar@{-}[d] \ar@{-}[dr] & \\
*+[F-,]{\pted{1}|\pted{2}3} \ar@{-}[drrr] &
*+[F-,]{\pted{1}|2\pted{3}} \ar@{-}[drr] &
*+[F-,]{\pted{2}|\pted{1}3} \ar@{-}[dr] & &
*+[F-,]{\pted{2}|1\pted{3}} \ar@{-}[dl] &
*+[F-,]{\pted{3}|\pted{1}2} \ar@{-}[dll] &
*+[F-,]{\pted{3}|1\pted{2}} \ar@{-}[dlll] \\
& & & *+[F-,]{\pted{1}|\pted{2}|\pted{3}} & && }$$
\caption{\label{figure1} The poset $\Pi^A_{3}$}
\end{figure}

The order relation on the set of pointed partitions of $[n]$ is
defined as follows. The underlying partitions must be related by
the refinement order of partitions and the set of pointed elements
of the finer partition $\pi$ must contain the set of pointed
elements of the other one $\nu$. In this case, one gets $\pi
\leqslant \nu$. For instance, one has $\{\pted{1}\}\{\pted{3}\}\{2
\pted{4} \}\leqslant \{\pted{1} 3\}\{2 \pted{4}\}$. We denote this
poset by $\Pi^A_{n}$. The example of $\Pi^A_{3}$ is given in
Figure~\ref{figure1}.

From the very definition of pointed partitions as non-empty sets
of pointed sets, one can see that the exponential generating
function for the graded cardinalities of the posets of pointed
partitions of type $A$ is given by $\frac{e^{x u e^u}-1}{x}$.
Indeed the exponential generating function for the pointed sets is
$u e^u$ and the exponential generating function for the non-empty
sets is $e^u-1$. The additional variable $x$ takes the number of
parts into account.

The symmetric group on $[n]$ acts by automorphisms on the poset
$\Pi^A_{n}$. The set of maximal elements has cardinality $n$ and
the action of the symmetric group is the natural transitive
permutation action. Hence all maximal intervals are isomorphic as
posets. We denote this poset by $\Pi^A_{n,\,1}$.
% see Figure~\ref{figure2} for an example.

The main property of this family of posets is the following one.

\begin{proposition}
  Each interval of $\Pi^A_{n}$ is isomorphic to a product of
  posets $\Pi^A_{\lambda_1,\, 1}\times \cdots \times
  \Pi^A_{\lambda_k,\, 1}$, where $\lambda_1+\cdots+\lambda_k \leqslant
  n$.
\end{proposition}
\begin{proof}
  First, any interval can be decomposed into a product according to
  the parts of the coarser partition. One can therefore assume that
  the maximal element of the interval is a single block. One can then
  replace, in each element of the interval, each block of the minimal
  element by a single element. This provides a isomorphism with some
  interval $\Pi^A_{\lambda,\,1}$.
\end{proof}

Hence the Möbius number of a pointed partition is the product of
the Möbius numbers of its parts. This property will allow us to
work by induction.

In general, the pointed partition posets are not lattices. They
are bounded below, pure posets and the rank of a pointed partition
$\pi$ of $[n]$ is equal to $n$ minus the number of blocks of
$\pi$.

\subsubsection{Type $B$}

Let us define the pointed partition poset of type $B_n$.\\

Let $[-n]$ be the set $\{-1,\dots,-n\}$. Recall first the
description of the usual partitions of type $B_n$. They are the
partitions of the set $[n]\sqcup [-n]$ such that there is at most
one block containing some opposite indices and the other blocks
come in opposite pairs. The block with opposite elements is called
the \emph{zero block}.\\

\begin{definition}[Pointed partition of type $B$]
\rm A \emph{pointed partition of type $B_n$} is a partition of
type $B_n$ together with the choice of an element of the zero
block and the choice for each pair of opposite blocks of a pair of
opposite elements. The chosen elements are called \emph{pointed}.
\end{definition}

For example, $\{\pted{3}, -2\}\{\pted{-1}, 1\}\{2, \pted{-3}\}$ is
a pointed partition of type $B_3$.

\begin{figure}[h]
$$\xymatrix@M=1pt@R=60pt@C=3pt{ & & *+[F-,]{\scriptstyle \pted{-2} {-1} {1} {2}} &
*+[F-,]{\scriptstyle{-2} \pted{-1} {1} {2}} & &
*+[F-,]{\scriptstyle{-2} {-1} \pted{1} {2}}
& *+[F-,]{\scriptstyle{-2} {-1} {1} \pted{2}} & &  \\
 *+[F-,]{\scriptstyle\pted{-1}|
\pted{-2} {2}|\pted{1}} \ar@{-}[drrrr] \ar@{-}[urr]  \ar@{-}[urrr]
\ar@{-}[urrrrr] & *+[F-,]{\scriptstyle\pted{-2}| \pted{-1}
{1}|\pted{2}} \ar@{-}[ur] \ar@{-}[urr] \ar@{-}[urrrrr]
\ar@{-}[drrr] & *+[F-,]{\scriptstyle\pted{-2} {-1} | {1} \pted{2}}
\ar@{-}[drr] \ar@{-}[u] \ar@{-}[urrrr] & *+[F-,]{\scriptstyle{-2}
\pted{-1} | \pted{1} {2}} \ar@{-}[dr] \ar@{-}[u] \ar@{-}[urr] &
&*+[F-,]{\scriptstyle{-2} \pted{1} | \pted{-1} {2}}\ar@{-}[dl]
\ar@{-}[u] \ar@{-}[ull] &*+[F-,]{\scriptstyle\pted{-2} {1} | {-1}
\pted{2}} \ar@{-}[dll] \ar@{-}[u] \ar@{-}[ullll] &
*+[F-,]{\scriptstyle\pted{-1}| {-2}
\pted{2}|\pted{1}}\ar@{-}[dlll] \ar@{-}[ul] \ar@{-}[ull]
\ar@{-}[ullll] &*+[F-,]{\scriptstyle\pted{-2}| {-1}
\pted{1}|\pted{2}} \ar@{-}[dllll] \ar@{-}[ull] \ar@{-}[ulll]
\ar@{-}[ullllll]
\\
 & & & & *+[F-,]{\scriptstyle\pted{-2}|\pted{-1}|\pted{1}|\pted{2}} & &&
}$$ \caption{\label{figure2} The poset $\Pi^B_{2}$}
\end{figure}
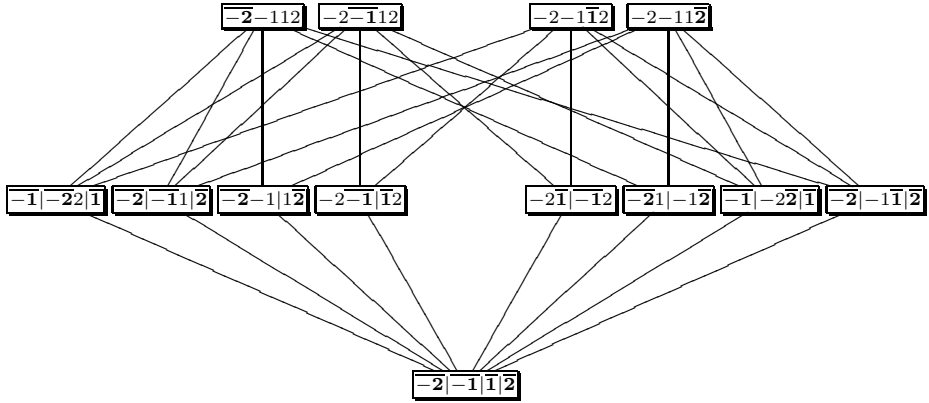

The order relation is as follows. The underlying partitions must
be related by the refinement order of partitions and the set of
pointed elements of the finer partition must contain the set of
pointed elements of the other one. For instance, one has
$\{\pted{-3}\}\{\pted{-2}\}\{\pted{-1}, 1\}\{\pted{2}\}\{
\pted{-3}\} \leqslant\{\pted{3}, -2\}\{\pted{-1}, 1\}\{2,
\pted{-3}\}$. We denote these posets by $\Pi^B_{n}$. The example
of $\Pi_2^{B}$ is given in Figure~\ref{figure2}.

The hyperoctahedral group of signed permutations of $[n]$ acts by
automorphisms on the poset $\Pi^B_{n}$. The set of maximal
elements has cardinality $2n$. Once again, the transitive action of
the hyperoctahedral groups on maximal intervals shows that they
are all isomorphic as posets. We denote this poset by
$\Pi^{B'}_{n}$.

In general, the pointed partition posets of type $B$ are not
lattices. They are bounded below, pure posets and the rank of a
pointed partition $\pi$ of type $B$ is equal to $n$ minus the
number of pairs of opposite blocks of $\pi$.

\begin{remark}
  \label{poset-beta}
  We consider also a variation of pointed partitions of type $B$ such
  that no element is pointed in the zero block. Such partitions are
  called \emph{pointed partitions of type $\beta$}. The partial order
  between pointed partitions of type $\beta$ is defined like the order
  in $\Pi^B_{n}$. We denote these graded posets by $\Pi^\beta_n$.
\end{remark}

%% \subsection{Type $D$}

%% PROBLEM HERE : DOES NOT WORK !

%% How to define a set of cardinality $2n-2$ acted upon transitively by
%% $D_n$ ?

%% Let us define the pointed partition poset of type $D_n$ ??

%% Let $[-n]$ be the set $[-1,\dots,-n]$. Recall first the description of
%% the partitions of type $D_n$. They are the partitions of type $B_n$
%% such that the zero block contains at least two pairs of opposite
%% indices.

%% A pointed partition of type $D_n$ is a partition of type $D_n$
%% together with the choice of an element of the zero block and the
%% choice for each pair of opposite blocks of a pair of opposite
%% elements.

%% The order relation as follows. The underlying partitions must be
%% related by the refinement order of partitions and the set of pointed
%% elements of the finer partition must contains the set of pointed
%% elements of the other one.

\subsection{Characteristic polynomials}

\subsubsection{Characteristic polynomials in type $A$}

The aim of this section is to compute the characteristic
polynomials of the posets of pointed partitions of type $A$. The
proof uses the subposets of pointed partitions for which a fixed
subset of $[n]$ is contained in the set of pointed elements. Up to
isomorphism, these subposets only depend on the cardinality of the
fixed subset of pointed indices. For $1 \leq i \leq n$, let
$\Pi^A_{n,i}$ be the poset where the indices in $[i]$ are pointed.
This poset has rank $n-i$. Recall that the pointed partition poset
of type $A_{n-1}$ is denoted by $\Pi^{A}_{n}$ and has rank $n-1$.

\begin{theorem}
  For $1 \leq i \leq n$, the characteristic polynomial of $\Pi^{A}_{n,i}$ is
  \begin{equation}
    \chi^{A}_{n,i}(x)=(x-i)(x-n)^{n-1-i},
  \end{equation}
  and its constant term $C^{A}_{n,i}$ is $(-1)^{n-i} i\, n^{n-1-i}$.
  The characteristic polynomial of $\Pi^A_{n}$ is
  \begin{equation}
    \chi^{A}_{n}(x)=(x-n)^{n-1}.
  \end{equation}
\end{theorem}
\begin{proof}
  The proof proceeds by induction on $n$. The statement of the Theorem
  is clearly true for $n=1$.

  From now on fix $n\geq 2$ and assume that the Theorem has been
  proved for smaller $n$. Then the proof is by decreasing induction on
  $i$ from $n$ to $1$. In the case $i=n$, the poset $\Pi^{A}_{n,n}$
  has just one element, so its characteristic polynomial is $1$, which
  is the expected value.

  So assume now that $i$ is smaller than $n$. Suppose first that $i$
  is at least $2$. Using the decomposition of a partition into its
  parts, which gives a product for the Möbius number, the constant
  term $C^{A}_{n,i}$ is given by
  \begin{equation}
    \label{recuC}
    \sum_{n_1,\dots,n_i} \prod_{j=1}^{i} \frac{C^{A}_{n_j,1}}{(n_j-1)!}(n-i)!,
  \end{equation}
  where the sum runs over integers $n_j \geq 1$ with sum $n$.

  Since $n_j < n$ for all $j$, we know by induction that
  $C^{A}_{n_j,1}=(-1)^{n_j-1}n_j^{n_j-2}$ and Lemma \ref{dominantA}
  then allows to compute the resulting sum:
  \begin{equation}
    C^{A}_{n,i}=(-1)^{n-i} i\, n^{n-1-i},
  \end{equation}
  as expected, when $i$ is at least $2$.

  Let us now compute $\chi^{A}_{n,i}$. By Möbius inversion on subsets
  of $[n]$ strictly containing $[i]$, it is given by
  \begin{equation}
    i\, n^{n-1-i}(-1)^{n-i}+\sum_{[i] \varsubsetneq S \subseteq [n]}
    (-1)^{|S|+i+1} \chi^{A}_{n,|S|} x^{|S|-i}.
  \end{equation}
  By induction on $i$, one gets
  \begin{equation}
    i n^{n-1-i}(-1)^{n-i}+\sum_{j=i+1}^{n}
    (-1)^{j+i+1} \binom{n-i}{j-i} (x-j)(x-n)^{n-1-j}x^{j-i} .
  \end{equation}
  Then using Lemma \ref{induction_step}, one finds the expected
  formula for $\chi^{A}_{n,i}$.

  Let us consider now the case $i=1$.

  Here we can not use induction to compute $C^{A}_{n,1}$, as Formula
  (\ref{recuC}) becomes trivial when $i=1$. Instead we use the fact
  that the poset $\Pi^{A}_{n,1}$ is bounded, hence $\chi^{A}_{n,1}$
  must vanish at $x=1$ and this property characterizes the constant
  term if the others coefficients are known. So let us guess what the
  constant term is and check later that the result vanish at $x=1$.

  Let us therefore compute $\chi^{A}_{n,1}$ as before, assuming that
  $C^{A}_{n,1}=(-1)^{n-1}n^{n-2}$. By Möbius inversion, it is given by
  \begin{equation}
    n^{n-2}(-1)^{n-1}+\sum_{[1] \varsubsetneq S \subseteq [n]}
    (-1)^{|S|} \chi^{A}_{n,|S|} x^{|S|-1}.
  \end{equation}
  By induction on $i$, one gets
  \begin{equation}
    n^{n-2}(-1)^{n-1}+\sum_{j=2}^{n}
    (-1)^{j} \binom{n-1}{j-1} (x-j)(x-n)^{n-1-j}x^{j-1} .
  \end{equation}
  Then using Lemma \ref{induction_step} again, one finds the expected
  formula for $\chi^{A}_{n,1}$. This formula vanishes at $x=1$, so the
  guess for the constant term was correct.

  Let us consider now the case of $\Pi^A_{n}$. In this case, one just
  has to use our knowledge of the other characteristic polynomials.
  Let us now compute $\chi^{A}_{n}$ as before. By Möbius inversion, it
  is given by
  \begin{equation}
    \sum_{\emptyset \varsubsetneq S \subseteq [n]}
    (-1)^{|S|+1} \chi^{A}_{n,|S|}x^{|S|-1}.
  \end{equation}
  By the previous results, one gets
  \begin{equation}
    \sum_{j=1}^{n}
    (-1)^{j+1} \binom{n}{j} (x-j)(x-n)^{n-1-j}x^{j-1} .
  \end{equation}
  Then using Lemma \ref{induction_step} one last time, one finds the
  expected formula for $\chi^{A}_{n}$.

  This concludes the inductive proof of the Theorem.
\end{proof}

\begin{lemma}
  \label{induction_step}
  For $0 \leq i \leq n$, one has the following equality:
  \begin{multline}
    i\, n^{n-1-i}(-1)^{n-i}+\sum_{j=1}^{n-i}(-1)^{j+1}
    \binom{n-i}{j}(x-(j+i))(x-n)^{n-1-i-j}x^j\\=(x-i)(x-n)^{n-1-i}.
  \end{multline}
\end{lemma}
\begin{proof}
  Introduce a new variable $y$ to get an homogeneous identity of degree
  $n-i$. Then replace $x$ by $1$ and $y$ by $(1-y)/n$. The resulting
  identity is easy to check.
\end{proof}

\begin{lemma}
  \label{dominantA}
  For $1 \leq i \leq n$, one has the following equation:
  \begin{equation}
    \sum_{n_1,\dots,n_i} \prod_{j=1}^{i} \frac{n_j^{n_j-1}}{n_j!}
    =i\, \frac{n^{n-1-i}}{(n-i)!},
  \end{equation}
  where the sum runs over integers $n_j \geq 1$ with sum $n$.
\end{lemma}
\begin{proof}
  Classical, see for example Proposition 2.5 in \cite{zvonkine}.
\end{proof}

\subsubsection{Characteristic polynomials in type $B$}

Let us compute the characteristic polynomials of the posets of
pointed partitions in type $B$. The proof uses the subposets where
$i$ and $-{i}$ are pointed for $i$ in a fixed subset of $[n]$. Up
to isomorphism, these subposets only depend on the cardinality of
the fixed subset of pointed pairs of indices. Let $\Pi^{B}_{n,i}$
be the poset where the indices in $[i]$ and $[-i]$ are pointed. By
convention, let $\Pi^B_{n,0}$ denote the pointed partition poset
$\Pi^{B}_{n}$. Recall that $\Pi^{B'}_n$ denotes a maximal interval
in $\Pi^{B}_{n}$.

\begin{theorem}
  For $0 \leq i \leq n$, the characteristic polynomial of $\Pi^{B}_{n,i}$ is
  \begin{equation}
    \chi^{B}_{n,i}(x)=(x-2n)^{n-i},
  \end{equation}
  with constant term $C^{B}_{n,i}=(-2n)^{n-i}$. The
  characteristic polynomial of $\Pi^{B'}_{n}$ is
  \begin{equation}
    \chi^{B'}_{n}(x)=(x-1)(x-2n)^{n-1},
  \end{equation}
  and its constant term is $C^{B'}_{n}=(-1)^{n}(2n)^{n-1}$.
\end{theorem}
\begin{proof}
  The proof proceeds by induction on $n$. The statement of the Theorem
  is clearly true for $n=1$.

  From now on fix $n\geq 2$ and assume that the Theorem has been
  proved for smaller $n$. Then the proof is by decreasing induction on
  $i$ from $n$ to $0$. In the case $i=n$, the poset $\Pi^{B}_{n,n}$
  has just one element, so its characteristic polynomial is $1$, which
  is the expected value.

  So assume now that $i$ is smaller than $n$. Suppose first that $i$
  is at least $1$. Using the decomposition of a partition of type $B$
  into its parts, which gives a product formula for the Möbius number,
  the constant term $C^{B}_{n,i}$ is given by
  \begin{equation}
    \sum_{n_1,\dots,n_i,m}
    \prod_{j=1}^{i}
    \frac{C^{A}_{n_j,1}}{(n_j-1)!}
    \frac{C^{B'}_{m}}{m!}
    2^{n-m-i} \,2m \, (n-i)!,
  \end{equation}
  where the sum runs over integers $n_j \geq 1$ and an integer $m\geq
  0$ with sum $n$.

  Using the results for type $A$, induction on $n$ to know
  $C^{B'}_{m}$ and Lemma \ref{dominantB} to compute the resulting sum,
  one gets that
  \begin{equation}
    C^{B}_{n,i}=(-2n)^{n-i},
  \end{equation}
  as expected.

  Let us now compute $\chi^{B}_{n,i}$. By Möbius inversion, one has
  the following equation:
  \begin{equation}
    \chi^{B}_{n,i}=(-2n)^{n-i}+\sum_{[i] \varsubsetneq S \subseteq
      [n]}(-1)^{|S|+i+1} \chi^{B}_{n,|S|} x^{|S|-i}.
  \end{equation}
  One gets by induction on $i$ that
  \begin{equation}
    \chi^{B}_{n,i}=(-2n)^{n-i}+\sum_{j=i+1}^{n}(-1)^{j+i+1}\binom{n-i}{j-i}
    (x-2n)^{n-j} x^{j-i},
  \end{equation}
  from which the expected formula follows through the binomial
  formula.

  There remains to compute $\chi^{B}_{n,0}$. For this, we need first
  to compute the characteristic polynomial $\chi^{B'}_{n}$ of a
  maximal interval $\Pi^{B'}_{n}$. Let us choose the maximal interval
  of elements where $n$ is pointed. Elements of this interval are of
  two distinct shapes: either $n$ is in the zero block or both $n$ and
  $-{n}$ are pointed. Let us split the computation of $\chi^{B'}_{n}$
  accordingly, as the sum of an unknown constant term, of
  $x \chi^{B}_{n,1}$ (already known) for the terms when $n$ and $-n$ are
  pointed and of the remaining terms when $n$ is in the zero block and
  the complement to the zero block is not empty.

  Let us compute this third part. There is a bijection between the set
  of such partitions and the set of triples $(S,\pi,\epsilon)$ where
  $S$ is a non-empty subset of $[n]\setminus\{n\}$, $\pi$ is a pointed
  partition of type $A$ on the set $S$ and $\epsilon$ is the choice of
  a sign for each element of $S$ up to complete change of sign of each
  block of $\pi$. The set $S$ is the positive half of the complement
  of the zero block and $\pi$ is the rest of the partition without its
  signs. Hence, one gets
  \begin{equation}
    \sum_{\emptyset \varsubsetneq S \subseteq [n-1]}
    \sum_{\pi \in \Pi^{A}_{|S|}} 2^{\rk(\pi)} \mu^A_{|S|}(\pi) C^{B'}_{n-|S|}
    x^{|S|-\rk(\pi)},
  \end{equation}
  where $\rk(\pi)$ is the rank in the poset $\Pi^{A}_{|S|}$.  This can
  be rewritten using induction on $n$ as
  \begin{equation}
    \sum_{j=1}^{n-1}
    \binom{n-1}{j}(-1)^{n-j}(2n-2j)^{n-j-1} 2^{j-1} x \chi^{A}_{j}(x/2).
  \end{equation}
  Then one can use the known results for type $A$ to obtain
  \begin{equation}
    \sum_{j=1}^{n-1}
    \binom{n-1}{j}(-1)^{n-j}(2n-2j)^{n-j-1} x (x-2j)^{j-1}.
  \end{equation}
  Using Lemma \ref{usefulB}, this is seen to be
  \begin{equation}
    -(x-2n)^{n-1}+(-2n)^{n-1}.
  \end{equation}

  As $\chi^{B'}_{n}$ has to vanish at $x=1$ because $\Pi^{B'}_{n}$ is
  bounded, one finds that its constant term is
  $(-1)^n(2n)^{n-1}$. Hence $\chi^{B'}_{n}$ is $(x-1)(x-2n)^{n-1}$ as
  expected.

  Now one can complete the proof by computing $\chi^{B}_{n,0}$ by
  Möbius inversion, just as for $\chi^{B}_{n,i}$ for $i>0$, because
  one now knows that $C^{B'}_{n}=(-1)^n (2n)^{n-1}$.
\end{proof}

\begin{lemma}
  \label{usefulB}
  For all $n \geq 1$, one has
  \begin{equation}
    \sum_{j=0}^{n-1}
    \binom{n-1}{j} (y+u j)^{n-j-1}(x-u j)^{j-1}=x^{-1}(x+y)^{n-1}.
  \end{equation}
\end{lemma}
\begin{proof}
  First set $m=n-1$ and rewrite it as
  \begin{equation}
    \sum_{j=0}^{m}\binom{m}{j}(y+j u)^{m-j}(x-j u)^{j-1}=x^{-1}(x+y)^{m}.
  \end{equation}
  Replace $y$ by $y+m$ and $u$ by $-1$. The identity becomes
  \begin{equation}
    \sum_{j=0}^{m}\binom{m}{j}(y+m-j)^{m-j}(x+j)^{j-1}=x^{-1} (x+y+m)^{m},
  \end{equation}
  which is one of many forms of the classical Abel binomial identity,
  see \cite{riordan} for example.
\end{proof}

\begin{lemma}
  \label{dominantB}
  For $1 \leq i \leq n$, one has the following equation:
  \begin{equation}
    \sum_{n_1,\dots,n_i,m} \prod_{j=1}^{i} \frac{n_j^{n_j-1}}{n_j !}
    \frac{m^m}{m!}
    =\frac{n^{n-i}}{(n-i)!},
  \end{equation}
  where the sum runs over integers $n_j \geq 1$ and an integer $m\geq
  0$ with sum $n$.
\end{lemma}
\begin{proof}
  Using the notations of Zvonkine \cite{zvonkine}, let
  \begin{equation}
    Y=\sum_{n \geq 1} \frac{n^{n-1}}{n!}u^n \text{   and   }
    Z=\sum_{n \geq 1} \frac{n^{n}}{n!}u^n.
  \end{equation}
  Then it is known (\cite[Prop. 2.5]{zvonkine}) that
  \begin{equation}
    Y^{i}= i \sum_{n \geq i} \frac{n^{n-i-1}}{(n-i)!}u^n.
  \end{equation}
  Applying the Euler operator $D=u\partial_u$, one gets
  \begin{equation}
    i \, Y^{i-1} Z= i \sum_{n \geq i} \frac{n^{n-i}}{(n-i)!}u^n.
  \end{equation}
  and the result follows because $Z=Y(1+Z)$.
\end{proof}

\subsubsection{Characteristic polynomials in type $\beta$}

Recall that $\Pi^{\beta}_n$ is the poset of partitions of type
$B_n$ where all blocks but the zero block are pointed, which was
defined in Remark \ref{poset-beta}. These posets appear as
intervals in the posets $\Pi^B$.

\begin{theorem}
  The characteristic polynomial of $\Pi^{\beta}_n$ is
  \begin{equation}
    \chi^{\beta}_{n}=(x-1)(x-(2n+1))^{n-1},
  \end{equation}
  and its constant term is $C^{\beta}_n=(-1)^n (2n+1)^{n-1}.$
\end{theorem}
\begin{proof}
  Let us prove the Theorem by induction on $n$. It is clearly true if
  $n=1$. Let us decompose the poset $\Pi^{\beta}_n$ according to the
  size $j$ of the complement of the zero block. Then the
  characteristic polynomial is
  \begin{equation}
    \sum_{j=0}^{n} \binom{n}{j} \sum_{\pi \in \Pi^A_j} \mu^A_j(\pi)
    2^{\rk(\pi)} C^{\beta}_{n-j} x^{j-\rk(\pi)},
  \end{equation}
  which can be rewritten as
  \begin{equation}
    \sum_{j=0}^{n} \binom{n}{j} 2^{j-1} \chi^{A}_j(x/2) C^{\beta}_{n-j} x.
  \end{equation}
  By known results on type $A$ and induction, the only unknown term is
  the constant term $C^{\beta}_{n}$, which is therefore fixed by
  the fact that the characteristic polynomial must vanish at $x=1$.
  Let us assume that this constant term has the expected
  value. One has to compute
  \begin{equation}
    \sum_{j=0}^{n} \binom{n}{j} x (x-2j)^{j-1} (-1)^{n-j} (2(n-j)+1)^{n-j-1} .
  \end{equation}
  Decomposing the binomial coefficient into
  $\binom{n-1}{j-1}+\binom{n-1}{j}$ and using twice the Abel binomial
  formula, one gets
  \begin{equation}
    (x-1)(x-(2n+1))^{n-1},
  \end{equation}
  which vanishes at $x=1$. This concludes the induction and the proof.
\end{proof}

\subsection{Homology}

The aim of this section is to compute the homology of the pointed
partition posets of type A and B. As a corollary, we get that the operad $\Pe$ is
Koszul over $\kk$.

For the different notions encountered in this section, we refer to
the article of A. Bj\"orner and M. Wachs \cite{BW}.

\subsubsection{Homology of $\Pi^A$}

Unlike the classical partition poset $\Pi_n$, which is a
semi-modular lattice, the pointed partition poset $\Pi^A_{n,\, 1}$
is not a lattice. Nevertheless, one has

\begin{lemma}
\label{thmAsemimod} For every $n\in \mathbb{N}^*$, the poset
$\Pi^A_{n,\, 1}$ is totally semi-modular.
\end{lemma}

\begin{proof}
  First, we prove that the poset $\Pi^A_{n,\, 1}$ is semi-modular for
  every $n\in \mathbb{N}^*$. Let $X$ and $Y$ be two different pointed
  partitions of $[n]$ covering a third pointed partition $T$. Denote
  the blocks of $T$ by $T=\{ T_1,\, T_2,\ldots,\, T_{k+1}\}$ and the
  pointed element of $T_i$ by ${t_i}$. Therefore, the pointed
  partitions $X$ and $Y$ are obtained from $T$ by the union of two
  blocks $T_i$ and $T_j$ and a choice of a pointed element between
  ${t_i}$ and ${t_j}$.  (We will often choose to denote these
  blocks by $T_1$ and $T_2$ for convenience). There are three possible
  cases.

  \begin{enumerate}

  \item The pointed partitions $X$ and $Y$ are obtained by the union
    of the same blocks $T_1$ and $T_2$. While ${t_1}$ is
    emphasized in $X$, ${t_2}$ is emphasized in $Y$. Since $X$ is
    different from $Y$ in the bounded poset $\Pi^A_{n,\, 1}$, $k$ must
    be greater than $2$. Consider the pointed partition $Z$ obtained
    from $T$ by the union of $T_1$, $T_2$ and $T_3$ where ${t_3}$
    is pointed. Therefore, $Z$ covers $X$ and $Y$.

  \item The pointed partition $X$ is obtained from $T$ by the union of
    $T_1$ and $T_2$ with ${t_1}$ emphasized and $Y$ is obtained by
    the union of $T_3$ and $T_4$ with ${t_3}$ emphasized.
    Consider the pointed partition $Z$ obtained from $T$ by the union
    of $T_1$ with $T_2$ and the union of $T_3$ with $T_4$ where
    ${t_1}$ and ${t_3}$ are emphasized. This pointed partition
    $Z$ covers both $X$ and $Y$.

  \item The pointed partition $X$ is obtained from $T$ by the union of
    $T_1$ and $T_2$ with ${t_i}$ emphasized ($i$=$1,\, 2$) and $Y$
    is obtained by the union of $T_2$ and $T_3$ with ${t_j}$
    emphasized ($j$=$2,\, 3$). We consider the pointed partition $Z$
    obtained by the union of $T_1$, $T_2$ and $T_3$. If $i$ is equal
    to $1$, we point out the element ${t_1}$ in $Z$. Otherwise, if
    $i$ is equal to $2$, we point out the element ${t_j}$ in $Z$.
    The resulting pointed partition $Z$ covers $X$ and $Y$.
  \end{enumerate}

  We can now prove that the poset $\Pi^A_{n,\,1}$ is totally
  semi-modular for every $n\in \mathbb{N}^*$. Let $[U,\, V]$ be an
  interval of $\Pi^A_{n,\,1}$. The poset $[U,\, V]$ is isomorphic to a
  product $\Pi^A_{\lambda_1,\, 1}\times \cdots \times
  \Pi^A_{\lambda_k,\, 1} $ of semi-modular posets. Therefore, $[U,\,
  V]$ is semi-modular.
\end{proof}

As a corollary, we get

\begin{theorem}
  The posets $\Pi^A_{n,\, 1}$ are CL-shellable and Cohen-Macaulay.
\end{theorem}
\begin{proof}
  CL-shellability follows from total semi-modularity by \cite[Corollary
  5.2]{BW}. Then the Cohen-Macaulay property follows from shellability.
\end{proof}

\begin{remark}
  We do not know whether the posets $\Pi^A_{n,\, 1}$ admit an
  EL-labelling.
\end{remark}

The following interesting relation to operads allows us to compute
the homology.

\begin{theorem}
  \label{ThmHomoA} The operad $\Pe$ is a Koszul operad over $\kk$ (the
  ring $\mathbb{Z}$ or any field). This is equivalent to the fact that
  the homology of the posets $\Pi^A_{n,\, 1}$ is concentrated in top
  dimension. Moreover, the homology of the posets $\Pi^A_{n}$ with
  coefficients in $\kk$ is given by the following isomorphism of 
  $\Sy_n$-modules 
  \begin{equation*}
    H_i(\Pi^A_{n})\cong \left\{
      \begin{array}{cc}
        \mathcal{RT}(n)^*\otimes \textrm{sgn}_{\Sy_n} & \textrm{if }
        i=n-1,\\
        0 & \textrm{otherwise},
      \end{array} \right.
  \end{equation*}
  where $\mathcal{RT}(n)$ is the $\Sy_n$-module induced by the free
  $\kk$-module on the set of \emph{rooted trees} (\emph{cf.}
  \cite{CL}).
\end{theorem}

\begin{proof}
  We use here the techniques described in \cite{homologie}. The set
  operad $\Pe$ gives rise to a family of posets $(\Pi_{\Pe}(n))_{n\in
    \mathbb{N}^*}$ which are isomorphic to $(\Pi^A_{n})_{n\in
    \mathbb{N}^*}$. And for every maximal pointed partition $\alpha$
  of the shape $\{1,\ldots,\, \pted{i},\ldots ,\, n\}$, the interval
  $[\zero,\, \alpha]$ is isomorphic to $\Pi^A_{n,\, 1}$.  Therefore,
  Theorem 10 of \cite{homologie} gives that the operad $\Pe$ is Koszul
  if and only if the posets $\Pi^A_{n,\, 1}$ are Cohen-Macaulay. Since
  these posets are CL-shellable, they are Cohen-Macaulay over the ring
  $\mathbb{Z}$ and over every field.

  Once again, the results of \cite{homologie} show that the homology
  groups of $\Pi^A_{n}$ of top dimension are isomorphic, as
  $\Sy_n$-modules, to the Koszul dual operad of $\Pe$, which is the
  operad $\Pli(n)$. And the operad $\Pli$ is known to be described by
  the representations of the symmetric groups on rooted trees.
\end{proof}

\subsubsection{Homology of $\Pi^B$}

Once again, we show that the posets $\Pi^{B'}_n$ are totally
semi-modular, CL-shellable and Cohen-Macaulay, for every $n\in
\NN^*$.\\

To do that, we need to understand the intervals of $\Pi^B_n$. Let
us introduce a variation of the posets $\Pi^B_n$ and $\Pi^\beta_n$
such that there is at most one pointed element in the zero block.
The partial order is defined like the order of $\Pi^B_n$ and
$\Pi^\beta_n$. The only difference is that if $\pi \leqslant \nu$
then the number of pointed elements in the zero block of $\nu$ is
greater than the number of pointed elements of the zero block of
$\pi$. We denote these posets by $\Pi^{\beta B}_n$. The maximal
interval of $\Pi^{\beta B}_n$ such that the maximal element is a
zero block with one pointed element is denoted $\Pi^{\beta B'}_n$.

\begin{proposition}
\label{produitPiB} Each interval of $\Pi^B_n$ is isomorphic to a
product of posets of the shape $\Pi^\beta_\lambda\times
\Pi^A_{\lambda_1,\, 1} \times \cdots \times \Pi^A_{\lambda_k,\,
1}$ or of the shape $\Pi^{\beta
  B'}_\lambda\times \Pi^A_{\lambda_1,\, 1} \times \cdots \times
\Pi^A_{\lambda_k,\, 1}$, where $\lambda+\lambda_1+\cdots
+\lambda_k \leq n$.
\end{proposition}

\begin{proof}
Let $[U,\, V]$ be an interval of $\Pi^B_n$. There are two possible
cases.

If the pointed elements of the zero blocks of $U$ and $V$ are the
same, then the refinement of the zero block of $V$ corresponds to
a poset of type $\Pi^\beta_\lambda$. The refinement of the other
blocks of $V$ with blocks of $U$ corresponds to posets of type
$\Pi^A_{\lambda_i,\, 1}$.

If the pointed elements of the zero blocks of $U$ and $V$ are
different, then the refinement of the zero block of $V$
corresponds to a poset of type $\Pi^{\beta B'}_\lambda$.
\end{proof}

\begin{lemma}
For every $n\in \NN^*$, the poset $\Pi^{B'}_n$ is totally
semi-modular.
\end{lemma}

\begin{proof}
With Proposition~\ref{produitPiB}, it is enough to show that the
posets $\Pi^\beta_n$ and $\Pi^{\beta B'}_n$ are semi-modular.\\

Let $X$ and $Y$ cover $T=\{-T_{k+1}, \ldots, -T_1,T_0,\,
T_1,\ldots, T_{k+1}\}$ in $\Pi^\beta_n$. If the zero block of $X$
and $Y$ is $T_0$, then the proof is the same as in the case
$\Pi^A_{n,1}$. Otherwise, there are three cases.

\begin{enumerate}

\item If $X$ is given by $-T_1 \cup T_0  \cup T_1$ and $Y$ by
$T_2\cup T_3$ (and $-T_2 \cup -T_3$), then we consider $Z$ defined
by $-T_1 \cup T_0 \cup T_1$ and $T_2\cup T_3$, with the same
choice of pointed elements.

\item If $X$ is given by $-T_1 \cup T_0  \cup T_1$ and $Y$ by
$T_1\cup T_2$ (and $-T_1\cup -T_2$), then we consider $Z$ defined
by $-T_2\cup -T_1 \cup T_0 \cup T_1 \cup T_2$, with no pointed
elements.

\item If $X$ is given by $-T_1 \cup T_0  \cup T_1$ and $Y$ by
$-T_2 \cup T_0  \cup T_2$, then we consider $Z$ defined by
 $-T_2\cup -T_1 \cup T_0 \cup T_1 \cup T_2$, with no pointed
elements.

\end{enumerate}

In each case, the partition $Z$ covers both $X$ and $Y$.\\

Let $X$ and $Y$ cover $T=\{-T_{k+1}, \ldots, -T_1,T_0,\,
T_1,\ldots, T_{k+1}\}$ in $\Pi^{\beta B'}_n$. The proof is mainly
the same as in the case $\Pi^\beta_n$, except some choices of
pointed elements in the zero blocks. The only new case is the
following one.

When $X$ and $Y$ are obtained from $T$ by $-T_1 \cup T_0 \cup T_1$
but a different choice of pointed element. Therefore, we consider
$Z$ defined by $-T_2\cup -T_1 \cup T_0 \cup T_1 \cup T_2$ with a
pointed element coming from $T_2$.
\end{proof}

Therefore, using the results of \cite{BW}, we have the following
theorem.

\begin{theorem}
The posets $\Pi^{B'}_n$ are CL-shellable and Cohen-Macaulay, for
$n\in \NN^*$.
\end{theorem}

\begin{remark}
Since the theory of operads is based on representations of the
symmetric groups $\Sy_n$, we can not use it here to compute the
homology groups of $\Pi^B_n$.
\end{remark}

\subsection{Extended pointed partition posets}

Let us define $\widehat{\Pi}^A_{n}$ as the bounded poset obtained
from $\Pi^A_{n}$ by adding of a maximal element $\one$.

\begin{theorem}
  The poset $\widehat{\Pi}^A_{n}$ is totally semi-modular,
  CL-shellable and Cohen-Macaulay. Its homology is concentrated in top
  dimension and has dimension $(n-1)^{n-1}$.
\end{theorem}
\begin{proof}
  Let us prove that this poset is semi-modular first.  The proof is
  essentially the same as for the poset $\Pi^A_{n}$. Only the first
  case can be different, when two blocks are gathered in two different
  ways and there is no other block. Then $\one$ covers both.

  Now any interval in $\widehat{\Pi}^A_{n}$ is either an interval
  in $\Pi^A_{n}$, hence semi-modular, or an interval $[\pi,\one]$.
  Such an interval is isomorphic to a poset
  $\widehat{\Pi}^A_{\lambda}$, hence semi-modular either.

  From this, one deduces the shellability and Cohen-Macaulay property.
  This implies the concentration of the homology in top dimension.

  The Möbius number of $\widehat{\Pi}^A_{n}$ is given by the
  opposite of the value at $x=1$ of the characteristic polynomial of
  $\Pi^A_{n}$. This gives the Euler characteristic, hence here the
  dimension of the homology.
\end{proof}

\begin{remark}
  The action of the symmetric groups on the top homology of the posets
  $\widehat{\Pi}^A_{n}$ certainly deserves further study. It should be
  related to the vertebrates (twice-pointed trees) and to the
  generators of the free pre-Lie algebras as Lie algebras.
\end{remark}

\subsection{Incidence Hopf algebra in type $A$}

Let us consider the set $\FF$ of isomorphism classes of all
intervals in all posets $\Pi^A_n$ for $n\geq 1$. Then it follows
from Section \ref{definitionPA} that a set of representatives of
isomorphism classes is provided by arbitrary (possibly empty)
products of the intervals $\Pi^{A}_{n,1}$ for $n\geq 2$. This
family of intervals is therefore closed under products and taking
subintervals. Such a family is called \emph{hereditary} in
\cite{schmitt}.

Hence one can introduce the incidence Hopf algebra of this family
of intervals. For short, let $a_n$ be the isomorphism class of
$\Pi^{A}_{n,1}$ for $n\geq 2$. Then a Hopf algebra structure is
defined on the polynomial algebra $H(\FF)$ in the $a_n$ by the
following coproduct:
\begin{equation}
  \Delta \,a_n:=\sum_{\pi \in \Pi^{A}_{n,1}} [\zero,\pi] \otimes [\pi,\one],
\end{equation}
where $[\,]$ denotes the isomorphism class of the underlying
interval.

By convention, let $a_1$ be the unit of the polynomial algebra in
the variables $a_n$ for $n \geq 2$. It corresponds to the class of
the trivial interval.

Let us decompose the coproduct according to the rank of $\pi$. One
gets

\begin{equation}
  \Delta \,a_n =\sum_{k=1}^{n} \left(\sum_{{\pi \in \Pi^{A}_{n,1}}\atop {\rk(\pi)=n-k}}
  \prod_{i=1}^{k} a_{\pi_i}  \right) \otimes a_k,
\end{equation}
where the $\pi_i$ denotes the size of the blocks of the pointed
partition $\pi$ ($\sum_{i=1}^k \pi_i =n)$.

Then decomposing the pointed partition $\pi$ into its parts, with
a distinguished part containing $1$, one gets
\begin{equation}
  \Delta \,a_n =\sum_{k=1}^{n} \left(
\sum_{\pi_1=1}^{n+1-k} \sum_{{\pi_2,\dots,\pi_k \geq
1}\atop{\pi_1+\dots+\pi_k=n}}
\frac{(n-1)!}{(\pi_1-1)!\pi_2!\dots\pi_k!} \frac{\prod_{i=2}^{k}
  \pi_i}{(k-1)!}
 \prod_{i=1}^{k} a_{\pi_i}
\right) \otimes a_k.
\end{equation}

This formula can be rewritten as
\begin{equation}
  \label{coproduit}
    \Delta \frac{a_n}{(n-1)!} =\sum_{k=1}^{n} \left(
\sum_{{\pi_1,\pi_2,\dots,\pi_k \geq 1}\atop{\pi_1+\dots+\pi_k=n}}
\frac{a_{\pi_1} \dots a_{\pi_k}}{(\pi_1-1)!\dots(\pi_k-1)!}
\right) \otimes \frac{a_k}{(k-1)!}.
\end{equation}

This has the following interpretation.

\begin{theorem}
  The incidence Hopf algebra of the family of pointed partition
  posets of type $A$ is isomorphic to the Hopf algebra structure on
  the polynomial algebra in the variables $(a_n)_{n \geq 2}$ given by
  the composition of formal power series of the following shape:
  \begin{equation}
  \label{formalpowerseries}
    x+\sum_{n \geq 2} a_n \frac{x^n}{(n-1)!}.
  \end{equation}
\end{theorem}
\begin{proof}
  This follows from the explicit formula (\ref{coproduit}) for the
  coproduct on the generators.
\end{proof}

As a corollary, the Möbius numbers of the intervals
$\Pi^{A}_{n,1}$ can be deduced from the fact that the inverse for
composition of $x \exp(x)$ is the Lambert $W$ function whose
Taylor expansion is known to be
\begin{equation}
  W(x)=\sum_{n \geq 1} (-1)^{n-1} n^{n-2} \frac{x^n}{(n-1)!}.
\end{equation}

% \begin{remark}
%   The posets of pointed partitions comes from the operad $\Pe$, as
%   explained in \cite{homologie}. P. van der Laan has described in
%   \cite{PvdL} a natural structure of Hopf algebra of the polynomial
%   algebra on $\bigoplus_{n\ge 2} (\Pe(n)_{\Sy_n})^*$, where the
%   diagonal comes from the operadic composition of $\Pe$. We can
%   identify the unique orbit of $\Pe(n)$ under the action of the
%   symmetric group $\Sy_n$ with $a_n$. This shows that the incidence
%   Hopf algebra $H(\FF)$ is isomorphic to the one introduced by P.  van
%   der Laan.

%   The formal power series $x \exp(x)$ is the Poincar\'e series of the
%   operad $\Pe$ and the series $W(x)$ is the Poincar\'e series of the
%   Koszul dual operad of $\Pe$, namely the operad $\Pli$.
% \end{remark}

\section{Multi-pointed partition posets}

In this section, we give the definition of the multi-pointed partition poset
of type A and its basic properties.

\subsection{Type $A$}

Let us define the multi-pointed partition poset of type $A_{n-1}$.

\begin{definition}[Multi-pointed partition]
A \emph{multi-pointed partition} of $[n]$ is a partition of $[n]$
together with the choice of a non-empty subset of each block,
called the \emph{pointed subset} of this block.
\end{definition}

\begin{figure}[h]
$$\xymatrix@R=80pt@C=15pt{ & *+[F-,]{\pted{1}23} &*+[F-,]{1\pted{2}3} &
*+[F-,]{12\pted{3}} & *+[F-,]{\pted{12}3} &
*+[F-,]{\pted{1}2\pted{3}}
& *+[F-,]{1\pted{23}} & *+[F-,]{\pted{123}} & \\
*+[F-,]{\pted{1}|\pted{2}3} \ar@{-}[ur] \ar@{-}[urr]
\ar@{-}[urrrr]  & *+[F-,]{\pted{1}|2\pted{3}} \ar@{-}[u]
\ar@{-}[urr] \ar@{-}[urrrr] & *+[F-,]{\pted{2}|\pted{1}3}
\ar@{-}[u] \ar@{-}[ul] \ar@{-}[urr] & *+[F-,]{\pted{2}|1\pted{3}}
\ar@{-}[ul] \ar@{-}[u] \ar@{-}[urrr]& *+[F-,]{\pted{3}|\pted{1}2}
\ar@{-}[ulll] \ar@{-}[ul] \ar@{-}[ur]&
*+[F-,]{\pted{3}|1\pted{2}}\ar@{-}[ulll] \ar@{-}[ull] \ar@{-}[ur]
& *+[F-,]{\pted{1}|\pted{23}} \ar@{-}[ulllll] \ar@{-}[u]
\ar@{-}[ur] & *+[F-,]{\pted{2}|\pted{13}} \ar@{-}[ulllll]
\ar@{-}[ull] \ar@{-}[u]& *+[F-,]{\pted{3}|\pted{12}}
\ar@{-}[ulllll] \ar@{-}[ullll]
\ar@{-}[ul] \\
& & & & *+[F-,]{\pted{1}|\pted{2}|\pted{3}} \ar@{-}[u] \ar@{-}[ul]
\ar@{-}[ull] \ar@{-}[ulll] \ar@{-}[ullll] \ar@{-}[ur] \ar@{-}[urr]
\ar@{-}[urrr] \ar@{-}[urrrr] & & & & }$$ \caption{\label{figure3}
The poset $\MPi^A_{3}$}
\end{figure}
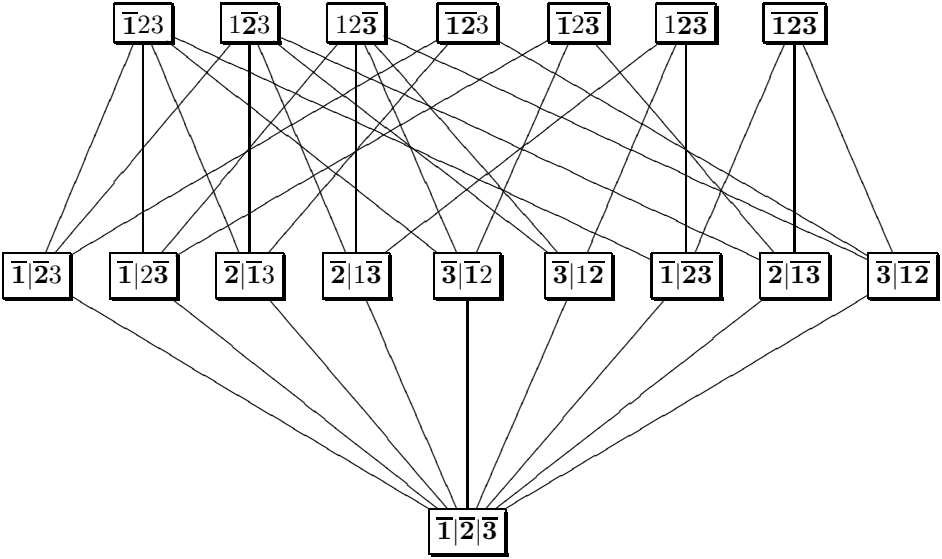

The order relation is as follows. First the underlying partitions
must be related by the refinement order of partitions. Then if two
partitions are related by the gathering of two blocks, the set of
pointed elements of the big block is either one of the sets of
pointed elements of the two small blocks or their union. For
instance, one has $\{\pted{1}2\}\{3 \pted{56}\}\{\pted{47}8\}
\leqslant \{123\pted{56}\}\{\pted{47}8\}$. The poset of
multi-pointed partitions of type $A_{n-1}$ is denoted by
$\MPi^A_{n}$. The example of the poset $\MPi^A_{3}$ is displayed
in Figure~\ref{figure3}.

As the multi-pointed partitions are just non-empty sets of pairs
made of a non-empty set and a set, the generating series for the
graded cardinality is given by $\frac{e^{x e^u(e^u-1)}-1}{x}$.

Of course, the symmetric group $\Sy_n$ acts on the poset
$\MPi^A_{n}$.

Let $\MPi^{A'}_{n,i}$ denote the maximal interval in $\MPi^A_{n}$
between $\zero$ and a multi-pointed partition with one block and
$i$ pointed elements. Clearly this does not depend on the choice
of the pointed elements.

The following proposition will play a crucial role in the sequel.

\begin{proposition}
\label{produitMPi}
  Each interval of $\MPi^A_{n}$ is isomorphic to a product of
  posets $\MPi^{A'}_{\lambda_1,\, \nu_1}\times \cdots \times
  \MPi^{A'}_{\lambda_k,\, \nu_k}$, where $\lambda_1+\cdots+\lambda_k
  \leqslant n$ and $1 \leqslant \nu_i \leqslant \lambda_i$.
\end{proposition}

\begin{proof}
  As for the pointed partition posets, any interval can be decomposed
  into a product according to the parts of the coarser partition.
  One can therefore assume that the maximal element of the interval is
  a single block. One can then replace, in each element of the
  interval, each block of the minimal element by a single element.
  This provides a isomorphism with some interval
  $\MPi^{A'}_{\lambda,\,\nu}$.
 \end{proof}

\subsection{Characteristic polynomials in type $A$}

Let us compute the characteristic polynomials of the posets of
multi-pointed partitions of type $A$. The proof uses the subposets
of elements where a fixed subset of $[n]$ is pointed. Up to
isomorphism, these subposets only depend on the cardinality of the
fixed subset of pointed indices. For $1 \leq i \leq n$, let
$\MPi^A_{n,i}$ be the poset where the indices in $[i]$ are
pointed. Let us denote by $\MPi^{A'}_{n,i}$ the maximal interval
under a partition with a single block and $i$ pointed elements. Up
to isomorphism, this does not depend on the choice of the pointed
elements. By convention, let $\MPi^{A'}_{n,0}$ denote the
multi-pointed partition poset $\MPi^A_{n}$ of type $A_{n-1}$.

Let us introduce the following convenient (if not traditional)
notation:
\begin{equation}
  \lbr n \rbr= \prod_{j=1}^{n} (x-j).
\end{equation}

Let us remark that $\MPi^{A}_{n,n}$ and $\MPi^{A'}_{n,n}$ are both
isomorphic to the classical partition poset of type $A_{n-1}$
whose characteristic polynomial is known to be $\lbr n-1 \rbr$.

\begin{theorem}
  For $1 \leq i \leq n$, the characteristic polynomial of $\MPi^{A}_{n,i}$ is
  \begin{equation}
    \Mchi^{A}_{n,i}(x)=\frac{x-2i}{x-i}
    \frac{\lbr i \rbr\lbr 2n-1 \rbr}{\lbr n+i \rbr},
  \end{equation}
  and its constant term $MC^{A}_{n,i}$ is $(-1)^{n-1} 2
  \frac{i!(2n-1)!}{(n+i)!}$. For $0 \leq i \leq n$, the characteristic
  polynomial of $\MPi^{A'}_{n,i}$ is
  \begin{equation}
    \Mchi^{A'}_{n,i}(x)=
    \frac{\lbr i \rbr\lbr 2n-i-1 \rbr}{\lbr n \rbr},
  \end{equation}
  and its constant term $MC^{A'}_{n,i}$ is $(-1)^{n-1}
  \frac{i!(2n-i-1)!}{n!}$.
\end{theorem}
\begin{proof}
  Let us prove the Theorem by recursion on $n$. It is clearly true for
  $n=1$. Let us now assume that it has been proved for smaller
  $n$. The proof goes in three steps.

  The first step is to compute $\Mchi^{A'}_{n,i}$ by decreasing
  recursion on $i$ for $i>0$. The statement is clear if $i=n$. Let us
  assume that the chosen pointed elements are $[i]$. The poset
  $\MPi^{A'}_{n,i}$ can be decomposed according to the size and number
  of pointed elements of the block $p_1$ containing $1$. Let $J$ be
  the intersection of $p_1$ with $[i]$. This is exactly the set of
  pointed elements of $p_1$. Let $S$ be the complement of $J$ in
  $p_1$, contained in $[n]\setminus [i]$. Then the result is
  \begin{equation}
    \sum_{[1]\subseteq J \subseteq [i]} \sum_{\emptyset \subseteq S
    \subseteq [n]\setminus[i]} \sum_{\pi \in
    \MPi^{A'}_{n-|S|-|J|,i-|J|}}
    \Mchi^{A'}_{|J|+|S|,|J|} \mu(\pi) x^{n-\rk(\pi)-|S|-|J|+1},
  \end{equation}
  where $\rk(\pi)$ is the rank in the poset $\MPi^{A'}_{n-|S|-|J|,i-|J|}$.
  Hence one gets the following equation for $\Mchi^{A'}_{n,i}$:
  \begin{equation}
    \sum_{j=1}^{i}\binom{i-1}{j-1}\sum_{s=0}^{n-i}\binom{n-i}{s}
    MC^{A'}_{j+s,j} \Mchi^{A'}_{n-j-s,i-j} x.
  \end{equation}
  The only unknown term is the constant term when $s=n-i$ and
  $j=i$. This coefficient is determined by the fact that
  $\Mchi^{A'}_{n,i}$ must vanish at $x=1$. So let us assume that it has
  the expected value and check later that the result vanish at
  $x=1$. One therefore has to compute
  \begin{multline}
    \sum_{j=1}^{i}\binom{i-1}{j-1}\sum_{s=0}^{n-i}\binom{n-i}{s}
    (-1)^{j+s-1}\frac{j!(2s+j-1)!}{(s+j)!}\\
    \frac{\lbr i-j \rbr \lbr 2n-2s-j-i-1\rbr}{\lbr n-j-s\rbr}
    x^{1}.
  \end{multline}

  Using Lemma \ref{convolution} to compute the inner summation on $s$
  and then Lemma \ref{facteur} to compute the remaining summation on
  $j$, one gets the expected formula for $\Mchi^{A'}_{n,i}$. As this
  formula vanish at $x=1$, the guess for the constant term was
  correct.

  The second step is to compute $\Mchi^{A}_{n,i}$ by decreasing
  recursion on $i$. By a decomposition as above according to the block
  containing $1$, one gets the following equation for
  $\Mchi^{A}_{n,i}$:
  \begin{equation}
    \sum_{j=1}^{i}\binom{i-1}{j-1}\sum_{s=0}^{n-i}\binom{n-i}{s}
    MC^{A}_{j+s,j} \Mchi^{A}_{n-j-s,i-j} x.
  \end{equation}
  The only unknown term is the constant term when $s=n-i$ and
  $j=i$. This coefficient is given by
  \begin{equation}
    \sum_{[i]\subseteq S \subseteq [n]} MC^{A}_{n,|S|}.
  \end{equation}
  This quantity is computed in Lemma \ref{Mdominant} and found to be
  as expected. To determine $\Mchi^{A}_{n,i}$, one therefore has to
  compute
  \begin{multline}
    \sum_{j=1}^{i}\binom{i-1}{j-1}\sum_{s=0}^{n-i}\binom{n-i}{s}
    (-1)^{j+s-1} 2 \frac{j!(2s+2j-1)!}{(s+2j)!}\\\frac{(x-2(i-j))}{x-(i-j)}
    \frac{\lbr i-j \rbr \lbr 2n-2s-2j-1\rbr}{\lbr n-2j-s+i\rbr}
    x.
  \end{multline}
  Using Lemma \ref{convolution} to compute the inner summation on $s$
  and then Lemma \ref{facteur} to compute the remaining summation on
  $j$, one gets the expected formula for $\Mchi^{A}_{n,i}$.

  The last step is to compute $\Mchi^{A'}_{n,0}$ by Möbius inversion on
  non-empty subsets of $[n]$. Indeed it is clear that
  \begin{equation}
    x \Mchi^{A'}_{n,0}=\sum_{\emptyset \varsubsetneq S \subseteq
    [n]}(-1)^{|S|+1} \Mchi^{A}_{n,|S|}.
  \end{equation}
  So we have to show that
  \begin{equation}
     x \frac{\lbr 2n-1 \rbr}{\lbr n \rbr} = \sum_{j=1}^{n} \binom{n}{j}
     (-1)^{j+1} (x-2j)\frac{\lbr
    j-1 \rbr \lbr 2n-1 \rbr}{\lbr n+j \rbr}.
  \end{equation}
  This can be restated as the vanishing of
  \begin{equation}
    \sum_{j=0}^{n} \binom{n}{j}(-1)^j (x-2j)\frac{\lbr
    j-1 \rbr}{\lbr n+j \rbr}.
  \end{equation}
  In hypergeometric terms, this is equivalent to
  \begin{equation}
    {_3}{F_2}\left(-n,y,y/2+1;y/2,y+n+1;1\right)=0.
  \end{equation}
  This follows from a known identity, see Appendix (III.9) in
  \cite{slater} for example, with $a=y$ and $b=y/2+1$.
  %Reference S3232 in Krattenthaler list.
\end{proof}

\begin{lemma}
  \label{convolution}
  For all $m\geq 0$ and $j,k \geq 1$, one has
  \begin{multline}
    \sum_{s=0}^{m}\binom{m}{s}\frac{(-1)^s j(2s+j-1)!}{(s+j)!}
    \frac{(x-k)\lbr 2(m-s)+k-1 \rbr}{ \lbr m-s+k
    \rbr}\\=
    \frac{(x-(j+k))\lbr 2m+j+k-1 \rbr}{ \lbr m+j+k
    \rbr}
  \end{multline}
\end{lemma}
% ok
\begin{proof}
  Once reformulated using the Pochhammer symbol, this is a direct
  consequence of the product formula associated to the following
  hypergeometric function:
  \begin{equation}
    \psi_x(\theta)={_2}{F_1}\left(\frac{x}{2},\frac{1+x}{2}; 1+x ;\theta\right)=
\left(\frac{2}{1+\sqrt{1-\theta}}\right)^{x},
  \end{equation}
  which can be found for example as Formula (5), page 101 in
  \cite{erdelyi}. More precisely, one takes the constant term with
  respect to $y$ in the Taylor coefficients with respect to $\theta$
  of the identity
  \begin{equation}
    \psi_y(\theta)\psi_x(\theta)=\psi_{x+y}(\theta).
  \end{equation}
\end{proof}

\begin{lemma}
  \label{facteur}
  For all $k \geq 1$, one has
  \begin{equation}
    x \sum_{j=1}^{k} \binom{k-1}{j-1}(-1)^{j-1}(j-1)!\lbr k-j-1 \rbr
    = \lbr k-1 \rbr.
  \end{equation}
\end{lemma}
% OK
\begin{proof}
  Once reformulated using Pochhammer symbols, this becomes equivalent to
  \begin{equation}
    {_2}{F_1}\left(-k,1;-y-k+1;1\right)=\frac{y+k}{k},
  \end{equation}
  which is just one instance of the Gauss identity.
\end{proof}

\begin{lemma}
  \label{Mdominant}
  One has the following identity:
  \begin{equation}
    \sum_{j=i}^{n} \binom{n-i}{j-i}\frac{j!(2n-j-1)!}{n!}=
    2 \frac{i!(2n-1)!}{(n+i)!}.
  \end{equation}
\end{lemma}
\begin{proof}
  Once reformulated using hypergeometric functions, this becomes
  equivalent to
  \begin{equation}
    {_2}{F_1}\left(i+1,-m;-2m-i+1;1\right)=2\frac{(m+i)!(2m+2i-1)!}
    {(m+2i)!(2m+i-1)!},
  \end{equation}
  which is just another instance of the Gauss identity.
\end{proof}

\subsection{Homology of $\MPi^A$}

In this section, we compute the homology of the posets
$\MPi^A_{n,\, 0}$. As a corollary, we get that the operad
$\Comtri$
is Koszul over $\kk$.\\

Once again, we show that each maximal interval of $\MPi^A_{n,\,
0}$ is totally semi-modular. Therefore, the homology of
$\MPi^A_{n,\, 0}$ is concentrated in top dimension. And we use
Koszul duality theory for operads to compute this homology in
terms of $\Sy_n$-modules.

\begin{lemma}
For every $n \in \NN^*$ and every $1\leqslant i \leqslant n$, the
poset $\MPi^{A'}_{n,\, i}$ is totally semi-modular.
\end{lemma}

\begin{proof}
Since each interval of $\MPi^{A'}_{n,\,i}$ is isomorphic to a
product $\MPi^{A'}_{\lambda_1,\, \nu_1}\times \cdots \times
\MPi^{A'}_{\lambda_k,\, \nu_k}$, where $\lambda_1+\cdots+\lambda_k
\leqslant n$ and $1 \leqslant \nu_j \leqslant \lambda_j$, it is
enough to show that every $\MPi^{A'}_{n,\, i}$ is a semi-modular
poset.

Let $X$ and $Y$ be two different multi-pointed partitions of $[n]$
covering a third multi-pointed partition $T$ in
$\MPi^{A'}_{n,\,i}$. Denote the blocks of $T$ by $T$=$\{
T_1,\ldots,\, T_{k+1}\}$ and the set of pointed elements of $T_i$
by $\mathcal{T}_i$. Therefore, the multi-pointed partitions $X$
and $Y$ are obtained from $T$ by the union of two blocks $T_i$ and
$T_j$ and a choice of a pointed elements between $\mathcal{T}_i$,
$\mathcal{T}_j$ or both. (We will often choose to denote these
blocks by $T_1$ and $T_2$ for convenience). There are three
possible cases.

\begin{enumerate}

\item The multi-pointed partitions $X$ and $Y$ are obtained by the
union of the same blocks $T_1$ and $T_2$. Since $X$ is different
from $Y$ in the bounded poset $\MPi^{A'}_{n,\, i}$, $k$ must be
greater than $2$. Consider the multi-pointed partition $Z$
obtained from $T$ by the union of $T_1$, $T_2$ and $T_3$ where the
set $\mathcal{T}_3$ is pointed. Therefore, $Z$ covers $X$ and $Y$.

\item The multi-pointed partition $X$ is obtained from $T$ by the
union of $T_1$ and $T_2$ with the set $\mathcal{X}_1$ of pointed
elements. The multi-pointed partition $Y$ is obtained by the union
of $T_3$ and $T_4$ with the set $\mathcal{Y}_2$ of pointed
elements. Consider the multi-pointed partition $Z$ obtained from
$T$ by the union of $T_1$ with $T_2$ and the union of $T_3$ with
$T_4$ where the set $\mathcal{X}_1 \cup \mathcal{Y}_2$ of element
is emphasized. This multi-pointed partition $Z$ covers both $X$
and $Y$.

\item The multi-pointed partition $X$ is obtained from $T$ by the
union of $T_1$ and $T_2$ and $Y$ is obtained by the union of $T_2$
and $T_3$ where $\mathcal{Y}_2$ denotes the set of pointed chosen
elements. If only the elements of $\mathcal{T}_1$ or the elements
of $\mathcal{T}_2$ are emphasized in $X$, then we built the same
kind of covering partition $Z$ as in the proof of
Lemma~\ref{thmAsemimod}. If the elements of $\mathcal{T}_1\cup
\mathcal{T}_2$ are pointed in $X$, we consider the multi-pointed
partition $Z$ given by the union $T_1 \cup T_2 \cup T_3$ where
only the elements of $\mathcal{T}_3$ are pointed, if the elements
of $\mathcal{T}_2$ are not pointed in $Y$, and where the elements
of $\mathcal{T}_1 \cup \mathcal{Y}_2$ are pointed otherwise. In
any case, the multi-pointed partition $Z$ covers $X$ and $Y$.
\end{enumerate}
\end{proof}

As a consequence, using results of \cite{BW}, we have
\begin{theorem}
  The posets $\MPi^{A'}_{n,\,i}$ are CL-shellable and Cohen-Macaulay.
\end{theorem}

Then the relation with operads allows us to determine the
homology, as follows.

\begin{theorem}
The operad $\Comtri$ of commutative trialgebras is a Koszul operad
over $\kk$ (the ring $\ZZ$ or any field). This is equivalent to the
fact that the homology of the posets $\MPi^{A'}_{n,\,i}$ is
concentrated in top dimension. Moreover, the homology of the
posets $\MPi^A_{n}$ with coefficients in $\kk$ is given by the
following isomorphism of $\Sy_n$-modules
$$H_i(\MPi^A_{n})\cong \left\{
\begin{array}{cc}
\Lie \circ \Mag (n)^*\otimes \textrm{sgn}_{\Sy_n} & \textrm{if }
i=n-1,\\
0 & \textrm{otherwise},
\end{array} \right.
$$
where $\Lie \circ \Mag (n)$ is the $\Sy_n$-module induced by
plethysm or equivalently by the operadic composition of the operad
$\Lie$, of Lie algebras, with the operad $\Mag$, of magmatic
algebras.
\end{theorem}

\begin{proof}
Once again, we use the Theorems proved in \cite{homologie}. The
partition posets associated to the operad $\Comtri$ are isomorphic
to the posets $\MPi^A_{n}$, for $n\in \NN^*$. Since the Koszul
dual operad of the operad $\Comtri$ is the operad $\Postlie$,
which is isomorphic as $\Sy$-module to the composition $\Lie \circ
\Mag$ (\emph{cf.} \cite{homologie}), we conclude by the same
arguments as in the proof of Theorem~\ref{ThmHomoA}.
\end{proof}

\bibliographystyle{plain}
\bibliography{pointed}

{\small \textsc{Institut Girard Desargues, Université Claude
Bernard Lyon 1, Bâtiment Braconnier, 21 avenue Claude
Bernard, 69622 Villeurbanne Cedex France}\\
E-mail address: \texttt{chapoton@igd.univ-lyon1.fr}\\
URL: \texttt{http://igd.univ-lyon1.fr/$\sim$chapoton}}\\

{\small \textsc{Laboratoire J.-A. Dieudonné, Universit\'e de Nice
Sophia-Antipolis, Parc Valrose, 06108 Nice Cedex, France}\\
E-mail address: \texttt{brunov@math.unice.fr}\\
URL: \texttt{http://math.unice.fr/$\sim$brunov}}

\end{document}